\documentclass[11pt]{amsart}

\usepackage{amsmath,amsthm, amssymb, color}
\usepackage{graphicx}
\usepackage{fullpage}
\usepackage[normalem]{ulem}

\usepackage{amsmath, amsthm, amssymb, bm, mathtools}
\usepackage{enumerate}
\usepackage[backend=biber, style=numeric-comp,maxnames=99, minnames=99]{biblatex}
\newcommand{\kk}{G}
\newcommand{\p}{\partial}

\newcommand {\eps}{\varepsilon}

\title{On the direct SLE-martingale approach to the integrability of Conformal Loop Ensembles}
\author{Jason Miller and Wendelin Werner}
\address{DPMMS, University of Cambridge}

\begin{document}

\begin {abstract}
This paper is about the derivation of ``integrability results'' (i.e., for example, explicit formulas for the probability of events involving a finite number of marked points) for Conformal Loop Ensembles, that are the (conjectured) scaling limits of discrete crititical planar lattice models from Statistical Physics.  We discuss the explicit solutions to the PDEs provided by the original SLE/CLE martingale approach that AI now makes it possible to uncover. This provides short and direct new proofs of a number of formulas that have recently been derived based on Liouville quantum gravity, conformal field theory or Coulomb gas methods, and also makes it possible to derive some new ones. We briefly illustrate this with the examples of the backbone exponent and annulus crossing probabilities for critical percolation, but this direct approach seems to apply to a wide class of observables of two-dimensional systems related to conformal field theories.
\end {abstract}

\maketitle

The main goal of the present paper is to point out the existence, in the present age of AI, of fairly direct routes to a number of celebrated formulas describing probabilities of events for conformal loop ensembles or, more generally, what are often known as ``correlation functions'' in the corresponding conformal field theories describing the (conjectural and in some cases proven) scaling limit of a number of critical lattice models (such as percolation, the Ising model, the random cluster models or $O(N)$-loop models). The idea is to look at the martingales obtained from the Markovian explorations of a conformal loop ensemble, and then to ``just'' solve the resulting PDEs explicitly. So this goes back to the original SLE martingales developed and used in \cite {LSW1,LSW2,LSW5},  applied to the explorations of conformal loop ensembles. The format of the paper is arguably a little unusual but the hope is to provide a clear global picture without technicalities.

\section {Historical narrative}

Describing the scaling limits of critical models from statistical physics has been the topic of numerous important works over the past four decades,  in both the theoretical physics and mathematics communities. The list of references could easily end up being more than ten pages long. One could, for instance, at least concatenate the bibliographies of the ICM proceedings contributions \cite {icmSch,icmSm,icmM,icmLG, icmSh,icmDDG,icmRV,icmHS}.
Let us provide here an (arguably biased and obviously simplified) historical account.

\subsection*{Pre-SLE predictions}
Ideas that first emerged in the physics community in the 1980s and 1990s include conformal field theory, quantum gravity (with the introduction of Polyakov's action and the Knizhnik--Polyakov--Zamolodchikov (KPZ) relation), and Coulomb gas methods (there are also some earlier important ideas, for instance the transfer matrix approach). Each of these ideas led to spectacular predictions about the behavior of the scaling limits of planar lattice models such as the critical Ising model, critical percolation, or self-avoiding walks (the general class of models relevant to our discussion here includes FK percolation, the Potts models, and the $O(N)$ models). In particular, the CFT approach connected to representation theory (the Virasoro algebra) and was able to produce closed formulas for the asymptotic behavior of probabilities of events involving a certain number of marked points (the ``correlation functions'' of the fields). The KPZ relation in Quantum Gravity and Coulomb gas methods provided the explicit values of critical exponents that typically describe the leading asymptotic order of probabilities in the discrete models. From the mathematical point of view, all these works were considered to be ``non-fully rigorous'' and treated as ``predictions,'' i.e., challenging open problems for mathematicians. Landmark papers include \cite {BPZ,KPZ,David,Nienhuis87,Cardy1,perco1} (see also the textbooks \cite {CFT1,Cardybook} and references therein).
While some of these works produced scaling limits of finite-dimensional marginals of the lattice models (such as correlations at several points), the idea that one might also want to look at the scaling limit of clusters or loops as random geometric and conformally invariant objects emerged in the 1990s; see, for instance, \cite {Aiz98,AB}.

\subsection*{SLE and SLE martingales}

The situation changed drastically in 1999 when Schramm \cite {Schramm} realized that it was possible to describe explicitly the (only possible conformally invariant) scaling limits of single interfaces in these statistical physics models using Loewner's equation with a Brownian driving function. The resulting explicit non-self-crossing fractal planar curves \cite{RS} were also amenable to some explicit computations. To describe the whole set of interfaces and not just one of them, one now uses a variant known as conformal loop ensembles (CLEs), which can be defined either via a branching tree of target-independent SLE variants \cite {ShCLE}, or via some sort of Poisson point process of SLE excursions (bubbles) obtained when discovering the loops one by one \cite {ShW}. In some sense, these CLEs completely define/describe the possible conformally invariant scaling limits of all the lattice models as random geometric objects. They in principle encapsulate all the information about these scaling limits.  Most of  the ``correlation functions'' appearing in the Physics literature could then in principle be interpreted as expectations of functions of this random configuration -- which makes them reminiscent of ``spatial finite-dimensional marginals''
[recall that prior to this CLE description, the ``full characterization'' of the scaling limit was in some sense actually done via the collection of these types of finite-dimensional marginals, which explains also why their derivation was so crucial at the time].

\subsection*{Example of percolation}

Let us illustrate some aspects of the previous narrative using the simplest possible example, namely critical percolation. For simplicity, we will consider (since this is the percolation model that was then proved to be conformally invariant by Smirnov \cite {Smirnov}; see also \cite {KSS} for a new derivation) critical site percolation on the triangular lattice. We emphasize, however, that the reference to a discrete lattice model is for presentation purposes only (they make it for instance  somewhat easier to describe the events that one is looking at), and that the content of the present paper is really about the study of the continuum objects (that arise in their scaling limit). In particular, we will not discuss here the question of the convergence (and the speed of convergence) of the discrete models to the continuum ones, which is of course an important and active field.

Building on the CFT approach,
Cardy famously predicted the limiting value of the probability of left-to-right crossings of rectangles \cite{Cardy1}; the values of the critical exponents describing the asymptotic behavior of the probability of existence of ``arms'' were predicted by several means: key players include den Nijs, Nienhuis, Duplantier, Saleur, Aizenman and others \cite{dN79,NRS80,perco1,perco2,perco3,perco4}. This included the prediction that the one-arm exponent is equal to $5/48$, i.e., that the probability that the origin is in a cluster of diameter greater than $R$ behaves like $R^{-5/48 + o (1)}$ as $R \to \infty$. The asymptotic behavior of the probability of events other than crossings of rectangles was also predicted, such as double crossings (top--bottom and left--right) of rectangles. This is known as Watts' formula \cite {Watts}.

The SLE approach provided a direct route to these results: as argued by Schramm already in his first paper \cite {Schramm}, based on crossing probabilities of squares (and this was then further confirmed by its locality property \cite {LSW1}), SLE$_6$ is the only possible conformally invariant scaling limit of percolation interfaces. The proof (by Smirnov \cite {Smirnov}) that percolation had a conformally invariant scaling limit reduced the problem to studying properties of SLE$_6$ directly, which turned out to be doable (in the case of the  rectangle crossings, one just needs to evaluate the probability that SLE$_6$ started from the bottom-left corner of the rectangle hits the top boundary of the rectangle before its right boundary), and therefore provided rigorous derivations of almost all the above predictions for discrete critical percolation and also made it possible to derive a couple of new ones.

At this point, in the very early 2000s, the situation was therefore the following: the scaling limit of critical percolation was fully described via the branching SLE$_6$ structure (which was well-defined thanks to the target-invariance property of SLE$_6$ derived in \cite{LSW1}; see also \cite{CN}), which essentially meant that probabilities of events involving a finite number of marked points could be described via relatively simple PDEs corresponding to the evolution of their conditional probabilities when exploring part of one interface by SLE$_6$ curves (for models other than percolation, this later also became possible with the description of conformal loop ensembles).

When these PDEs in fact involved only one (conformal) variable (and in some few cases for two variables), it was possible to immediately guess/find their solutions (and therefore get explicit expressions for the corresponding probabilities and exponents). This led to the derivation of arm exponents for percolation \cite {LSW5,SW}, Watts' formula \cite {Dub,ShWi}, and Schramm's percolation formula \cite {SchrammFormula}. The non-monochromatic $k$-arm exponent (for $k \ge 2$) was in this way confirmed/proved to be
$(k^2 - 1)/12$ (see \cite {SW}, building on \cite {LSW2}), and the one-arm exponent to be $5/48$~\cite{LSW5}.
The expression for these exponents were of a rather simple form (and in particular involved no trigonometric functions). As we will discuss below, things looked somewhat messier when the number of variables was higher (which loosely corresponded, in some sense, to events involving more points, i.e., events involving more  conformally invariant quantities).

One remark that is relevant to the present discussion is that in these derivations of the values of the critical exponents (such as arm exponents) as eigenvalues of second-order differential operators, one was in fact also obtaining the corresponding eigenfunction, which was first done in a somewhat ad hoc manner. This was easy in those cases because the eigenfunctions were essentially the simplest possible functions (like powers of the functions $\sin(\theta)$ or $x(1-x)$; see the derivations in \cite {LSW1,LSW2,LSW5}, for instance).

\subsection*{What was left unsolved/unfinished}
Let us, however, describe two examples in the percolation setting where the limiting expression of the probability of an event was described in terms of the unique solution to an explicit PDE, but where the solution to the PDE could not be written out explicitly:

\begin {itemize}
 \item The backbone exponent: this is the value $\xi$ such that the probability that there exist two disjoint open paths (for critical percolation) that join the origin to the circle of radius $R$ decays like $R^{-\xi + o(1)}$ (disjoint means that the paths use different sites apart from the origin). It was explained in Appendix B of \cite {LSW5} that this value was the leading eigenvalue of some SLE-related second-order differential operator in two variables, but the explicit value was not determined there. It was later shown in \cite{BN} that $\xi$ lies strictly between the polychromatic two-arm and three-arm exponents, i.e., $1/4 < \xi < 2/3$.

 \item Crossings of annuli: Cardy \cite {Cardy06} proposed in 2006 (following his earlier work \cite{Cardy2} from 2002) an explicit formula for the limit $F(r)$ (when the mesh size of the lattice goes to $0$) of the crossing probability of the annulus between two concentric circles of radii $r=e^{-m}$ and $1$ (or conformal annuli with the same modulus). The expression he obtained is very explicit and involves ratios of values of the Dedekind eta function. If one writes $E(m) = \eta (i m/ \pi) :=e^{-m/12}\prod_{n=1}^{\infty}(1-e^{-2nm}) = r^{1/12} \prod_{n=1}^\infty (1-r^{2n})$, then
 \begin {equation}
  \label {Cardy}
  F(r) = \sqrt{\frac32}\,
       \frac{E(3m)E(3m/4)}{E(m) E( 3m/2)}.
       \end {equation}
       Along the same lines as in \cite {LSW5},
       Dub\'edat \cite{Du} obtained a PDE approach/characterization of annulus crossing probabilities (in particular of the probability that there exist both an open and a closed crossing, for which Cardy had proposed a formula in \cite{Cardy2}), and the same approach also yields a PDE characterization of $F(r)$, but the connection to Cardy's formulas was left open.
\end {itemize}

More specifically,
 what the CLE approach provides is a PDE satisfied by the function describing the ``state'' during the exploration of a boundary-touching loop. So, for the case of the backbone exponent, this reduces the problem to determining the leading eigenfunction (i.e., the strictly positive ones) $G(\alpha, \beta)$ on the state space of pairs of angles $\alpha, \beta$ in $(0, 2 \pi)$ with $\alpha + \beta < 2 \pi$ (so the domain for the set of pairs of angles $(\alpha, \beta)$ is in fact a triangle) for some simple second-order operator with some rather simple (but not always so simple to handle) boundary conditions. The corresponding eigenvalue (up to sign) is then the backbone exponent. We will rewrite this PDE and these conditions in Section \ref {S21}.
In Appendix B of \cite {LSW5}, it is written:
\begin{quote}
\emph{So far, we have not been able to solve explicitely the PDE problem and to give an
explicit formula for the value of the exponent, but this might perhaps be doable by someone
more proficient in this field.}
\end{quote}
\noindent However, this approach was not pushed much further at that stage. One reason for this probably was that it was not clear (and many of us did actually have doubts about this) at the time that there would be an explicit formula at all for this exponent (there was no conjecture for it in the physics literature, nor any argument indicating that it should take a special value that can be expressed by a formula; and  simulations did not seem to point to a particular noticeable value either).

Similarly, using the SLE$_6$ exploration as in \cite {Du}, one sees that the annulus crossing probability could also be characterized via some rather simple PDE with two variables (one of which being the modulus of the annulus that remains to be explored, the other being the (conformal) angle $\theta$ of the closed portion of the explored boundary). More precisely, the probability $F(r, \theta)$ that the inner circle in the annulus is connected to a given arc of length $2 \pi - \theta$ on the outer unit circle is the solution to some explicit PDE with simple boundary conditions (we will recall this PDE and these boundary conditions in Section \ref {S22}). The annulus crossing probability is then $F(r)= F(r,0)$. Again, it seemed unclear how to proceed in order to derive the formula (\ref {Cardy}) for $F(r,0)$ (i.e., how to obtain the general form of $F(r,\theta)$ needed to then identify $F(r,0)$ along such lines). Here too, a more ``proficient'' and persistent person might have been able to succeed.

It is worth emphasizing that obtaining these descriptions of probabilities and exponents in terms of PDEs is {\em short and direct}. One just needs to implement the very same arguments that led to the identification of the critical arm exponents in \cite {LSW1,LSW2,LSW5}, with just an occasional little twist to determine the right boundary conditions (see the example of the two monochromatic arms that we will briefly detail in Section \ref{S21}). Also, while \cite {LSW5} dealt with SLE$_6$, i.e., the scaling limit of critical percolation, the strategy can easily be extended to the other conformal loop ensembles (which were later defined and studied in \cite{ShCLE,ShW}), replacing the radial SLE$_6$ by the target-independent SLE$_\kappa (\kappa-6)$ explorations~\cite{SchW05,SSW09}.

\subsection*{Some further developments up to 2020}

Meanwhile, the mathematical theory of conformally invariant random systems was enriched by success stories in a number of very substantial ways. Two main ones directly relevant to the present discussion (as they were producing formulas and results for functionals of CLEs or CFT correlation functions) are the following (there is also the important approach via discrete planar maps and random metrics, which we will not discuss here -- see \cite {icmLG,icmDDG} and the references therein; similarly, the relation to the Gaussian free field as initiated in \cite {SchSh2,SchSh1,Dub09,MS1} and the Brownian loop soup \cite {LW04,ShW} and restriction properties \cite {LSWr} are also not so directly relevant to the present discussion):

\begin {itemize}
 \item Building on the original vision of Sheffield's quantum zipper paper \cite {ShQZ}, Duplantier, Miller, and Sheffield \cite {DMS} explained how to make sense of the Liouville quantum gravity approach and the KPZ relation (see also \cite{DS11}), and provided a mathematically clear framework to understand a number of the LQG derivations and arguments originating in the physics literature. This gave rise to a number of results (also relating to random planar maps). The almost sure KPZ formula \cite {GHM} could then, for instance, be used to derive values of critical exponents (for instance, it gives new proofs of the double and cut point dimensions of SLE first computed in \cite {MiWu}, which are related to four-arm exponents) for CLEs in the Euclidean setting, in a way reminiscent of the original KPZ approach by Duplantier, Saleur, and others~\cite{Dup98}.
 \item The ``integrable'' Liouville approach, initiated by Rhodes, Vargas, David, and Kupiainen (see \cite {DKRV,KRV} and a complete reference list in \cite {icmRV}), which derived some of the fundamental CFT structural formulas (like the DOZZ formula) via suitable studies of integrals involving the Gaussian free field. This approach did not involve any SLE interfaces or decompositions, but worked at the level of correlation function identities.
\end {itemize}

Our papers \cite {MSW1,MSW2} on CLE on LQG were the CLE counterparts of the mating-of-trees results for SLE on LQG. There, we relate the explorations of these CLEs to branching-type structures of weighted trees that are amenable to some explicit computations (and there is actually an ``integrability'' literature on those trees).
In these papers, a type of critical exponents that was new to us showed up for the first time. Some of these exponents can be described in terms of colored percolation models (if one takes the critical percolation case), known as divide-and-color models or fuzzy Potts models~\cite{MV95,Hag01}. Roughly speaking, one colors each percolation cluster independently in red or blue with probabilities $p$ and $1-p$, and one then looks at the ``clusters of clusters'' that one obtains. Among the results derived there (this was done directly in the continuum setting in \cite {MSW1,MSW2} and then shown to apply to discrete colored percolation clusters in \cite {KL}) was the value of the half-plane exponent $u(p)$ describing the decay of the probability that, for critical bond percolation in the upper half-plane, there exists a ``monochromatic chain of clusters'' (i.e., a path of red sites) from the origin to the semicircle of radius $R$. Assuming that critical bond percolation is conformally invariant in the scaling limit, this probability decays like $R^{-u(p)+o(1)}$ as $R \to \infty$, where
\begin {equation}
\label {u(p)}
 u(p) =\frac{1}{6} \Bigl(2 - \frac 3 \pi \arctan \frac{p\sqrt{3}}{2-p} \Bigr) \Bigl(1 - \frac 3 \pi \arctan \frac{p\sqrt{3}}{2-p} \Bigr).
 \end {equation}
This type of formula for a percolation-type critical exponent appeared quite new and therefore puzzling. In the introduction of \cite {MSW2}, it was written:
\begin{quote}
\emph{The appearance of the $\arctan(\cdot)$ function in the above formula suggests that its derivation will
involve arguments somewhat different from those used to derive the ``usual'' critical exponents,
which are computed only using SLE martingales. Indeed, as we shall explain in the present
paper, these $\arctan$ type formulas appear to be a by-product of the decomposition of LQG-type
surfaces in terms of L\'evy trees.}
\end{quote}

At that time, it already appeared that these results were maybe not that far from questions related to the backbone exponent (the second author here remembers his then-student Matthis Lehmkuehler being quite optimistic about it).

\subsection*{Recent developments}

In a spectacular development,  Nolin, Qian, Sun, and Zhuang \cite {NQSZ,NQSZ2} managed to derive the value of the backbone exponent in 2023, more than three years ago now. The value turns out to be the unique solution $\xi$ in $(1/4,2/3)$ of the equation
\begin {equation}
 \label{backbone}
\sin\left(\frac{2\pi}{3}\sqrt{12 \xi + 1}\right) + \frac{\sqrt{3}}{4}\sqrt{12\xi + 1} = 0
\end {equation}
(numerically, $\xi \approx 0.35667$).
This came as a surprise to a large part of the community, as many of us did not expect that there would be an explicit formula for this exponent, and the form of the exponent was strikingly different from the polychromatic arm exponents. Their proof combines important inputs from all three aforementioned approaches (Euclidean CLE theory, CLE on LQG/mating of trees, Liouville CFT -- the Liouville integrability plays apparently an essential role), and builds on a number of related contributions (including papers involving Ang or Holden; see, e.g., \cite{AHS,ARS2}); it is a remarkable achievement to have found such a path to this result. Given all the input that feeds into the proof, a fully self-contained write-up of the proof along these lines would be several hundred pages long.

Building on the same ideas and tools, a number of explicit formulas for probabilities involving CLEs (and therefore critical percolation scaling limits for the CLE$_6$ case) were derived. This includes new formulas (which had also not been explicitly stated in the physics literature before and which came as a surprise to most of the community), as well as the aforementioned annulus crossing formula of Cardy \cite{AnnulusProof}. This group of results has been dubbed ``integrability of CLE'' and one now indeed has many new explicit formulas for probabilities and exponents involving CLEs; see \cite {CLEintegrability,integ2,Misc1,Misc2,Misc3,int4}.

A number of further explicit formulas for observables that can be viewed as CLE observables were also recently derived by other means, related to the GFF coupling or Coulomb gas ideas (see, for instance, \cite {AL,J1,J2,J3}).

\subsection*{Revisiting the SLE martingale approach}
These recent developments raise very naturally the idea to revisit the question of the existence of the aforementioned simple SLE-PDE-only derivations of these formulas. Indeed, feeding the backbone result of \cite {NQSZ} into the original SLE-based description shows that this particular exponent is the leading eigenvalue of the corresponding second-order differential operator. It therefore seemed natural to immediately suspect (and we were surely not the only ones) that if the eigenvalue is that explicit, then there must also exist an explicit expression for the eigenfunction $G$. And if one could guess this explicit expression for the (non-negative) eigenfunction, then to determine the backbone exponent, it would suffice to check that this particular non-negative expression is an eigenfunction (with the right boundary conditions, including the ``balance'' condition at the corner $(0,0)$ discussed in Section~\ref{S21}), as this then implies by Perron--Frobenius considerations that its eigenvalue is indeed the leading eigenvalue of the differential operator and therefore equal to the backbone exponent; this would provide a direct derivation of the value of the backbone exponent. In other words, knowing that there actually is an explicit expression for the backbone exponent strongly suggests that ``{\em someone more proficient}'' in handling solutions to PDEs could actually provide a short SLE-based proof.

Similarly, that {\sl proficient person} would be able to derive the explicit expression for $F(r, \theta)$ in the annulus crossing formula from the PDE. Again, for a complete proof, one would just have to check that this expression is the solution to the PDE (obtained along the lines of \cite {Du}) that characterizes it in order to show that it is indeed equal to $F$, and therefore derive the annulus crossing formula proposed by Cardy and proved in \cite {AnnulusProof}.

Similar alternative SLE-based proofs of the various integrability results might also be around the corner. Note that, very naturally, this approach leads to hierarchies of multi-point functions (one explores a domain with $N$ marked points until the exploration splits them into $j$ and $N-j$ points) and is naturally suited to fields with local scaling exponents as well.

To our knowledge,  these PDEs had not been solved explicitly. But given the type of task required here (i.e., finding explicit expressions for solutions to second-order PDEs), it looks like something that one can reasonably ask AI to do, and at which AI is likely to be successful as it would know where to find the adequate tricks. Recall that one only needs to come up with an educated guess of what the explicit eigenfunction is. Then, in order to obtain a full proof, one just needs to check that this particular function solves the PDE (with the appropriate boundary conditions).

\medbreak

So, as a proof of concept, we will present what happens when prompting ChatGPT to find the solutions to the PDE problems described in \cite {LSW5} and obtained along the lines of \cite {Du}, in order to get ``direct'' proofs of the corresponding results as explained above. Sure enough, the explicit expressions for $G$ and $F$ were produced after a couple of minutes, and in fact, ChatGPT was also able to explain how to derive them, showing that this would actually have been accessible to a suitably motivated human.

The next section, in which we provide these explicit expressions, is meant to point out that proofs of this kind are now around. We will not write out the detailed proofs of these results here (we note however that AI did indeed produce Lean versions), as our point is really to make the community at large aware of the fact that these direct CLE-exploration proofs of these integrability results do work and are now available, rather than to start an extensive list of all the formulas (including new ones; note that the formulas for $F(r, \theta)$ and $G$ seem to be new) that one can obtain in this way.

The present narrative highlights the case of critical percolation (and therefore CLE$_6$) for presentation purposes only, but this general strategy applies to CLEs in general; one just replaces the SLE$_6$ explorations by an SLE$_\kappa (\kappa -6)$ exploration of CLEs. In particular, the examples that we described in this section also have versions/variants for CLE$_\kappa$ explorations for all $\kappa \in (4,8)$. One can expect a number of similar derivations for CLE$_\kappa$ for $\kappa \in (8/3, 4]$ as well, and for the BCLEs (boundary conformal loop ensembles) introduced in \cite {MSW}. Mind that a number of the aforementioned recent results about integrability of CLE were in fact also derived for general CLEs -- the backbone exponent  is for instance derived in \cite {NQSZ} for all $\kappa \in (4,8)$, not just for $\kappa=6$.

\section {Two examples}

\subsection {The backbone eigenfunction}
\label {S21}
We give here some details about the percolation backbone exponent derivation.

\subsubsection* {The PDE problem obtained in \cite {LSW5} via SLE exploration}

We first quote from Appendix B of \cite {LSW5}. Note that while the derivation leading to this is somewhat informal, the arguments can be made rigorous using the same ideas as in that paper (we do not quite quote verbatim, as we here correct a minus sign typo in \cite {LSW5}, and we write $G$ for the function denoted by $G_2$ there):

\begin{quote}
\emph{(...) one looks for the non-negative eigenfunctions of the more symmetric operator
$$
3 (\partial_\alpha^2 - 2 \partial_\alpha \partial_\beta
+ \partial_\beta^2 )
+ \cot\bigl(\frac \alpha 2\bigr) \p_\alpha
+ \cot\bigl(\frac \beta 2\bigr) \p_\beta
$$
 with boundary conditions (for $\alpha , \beta>0$, $\alpha+\beta<2\pi$),
$$
\kk (\alpha, 2 \pi - \alpha) = 0 , \
\kk (0, \beta) = 1
\hbox { and }
\partial_\beta \kk ( \alpha, 0) = 0 .$$}
\end{quote}

(In particular, $G$ vanishes at the corner $(2\pi,0)$, while at the corner $(0,2\pi)$, where the constant values $1$ and $0$ on the two sides meet, $G$ is not continuous -- so some more explanations are in fact required to say what one means here, see the discussion below. The choice of the value $1$ of $G(0, \beta)$ is not important, as any multiple of $G$ is also an eigenfunction -- choosing this value to be $1$ is just a way to identify $G$ uniquely).

These boundary conditions do not actually specify the eigenvalue uniquely, as they do not exclude the introduction of a ``killing/birth'' factor at the $(\alpha, \beta)=(0,0)$ corner which, in the case where this factor is positive, would kill the exploration at a rate proportional to the local time at $0$ of $\alpha$. It is, however, not difficult to introduce the additional ``balance'' condition that then identifies the eigenfunction and eigenvalue uniquely.

Let us give a few details about where this PDE comes from. The basic idea is very simple: one launches a radial SLE$_6$ process $(\gamma_t, t \ge 0)$ in the unit disk $D$ toward the origin, parametrized by minus the logarithm of the conformal radius of the connected component $D_t$ of $D \setminus \gamma [0,t]$ containing the origin (and the conformal radius is comparable up to multiplicative constants to the distance between $\gamma[0,t]$ and the origin). This process explores pieces of the CLE$_6$ loops. Let $\alpha_t$ and $\beta_t$ be the (conformal) angles (as seen from the origin) of the two boundaries of the {\em currently explored CLE loop corresponding to a {\em closed} cluster}, where $\alpha_t$ corresponds to its closed side and $\beta_t$ to its open side.

When both $\alpha_t, \beta_t$ are in $(0, 2 \pi)$ (note that $\alpha_t + \beta_t \le 2 \pi$), the process evolves according to the radial SLE$_6$ evolution. When $\beta$ hits $0$, the process is just reflected (this is just as in the one-arm estimate), but the definition means that whenever $\alpha_t$ hits $0$, $\beta$ is reset to $0$ (because one starts exploring a new CLE loop, so the process restarts at the corner $(0,0)$ of the triangle). These boundary conditions completely describe the process $(\alpha_t, \beta_t)$.

At the time $T$ at which $\alpha_t$ or $\beta_t$ hits $2 \pi$ (and this can happen without $\beta$ or $\alpha$ hitting $0$ just before $T$, so that the process will actually approach the corners of the triangle in a particular way from the inside of the triangle), it is easy to see that this will prevent any point in $D_T$ from being in $Q$ (loosely speaking defined to be the set of points for which a ``backbone event'' occurs): if $\alpha_T = 2 \pi$, then no open arm exists at all, as one has discovered a closed loop, and when $\beta_T = 2 \pi$, there is room for only one arm to go through the created pinch point. On the other hand, when $T$ has not happened yet, it is still ``possible'' for a backbone event to occur: one just needs two disjoint open paths, one reaching the open part of the boundary (corresponding to a conformal angle $2 \pi - \alpha_t$) and one reaching the part of the boundary corresponding to the angle $2 \pi - \alpha_t - \beta_t$ (this is because only one of the paths can then continue along the portion corresponding to $\beta_t$, except that whenever $\alpha$ hits zero, both directions on the $\beta_t$ arc become available, and one can indeed reset $\beta$ to $0$). So one indeed has to evaluate the asymptotic behavior of $P[T > t ]$ as $t \to \infty$.
This amounts to determining the leading eigenvalue of the operator associated with this killed Markov process (the probability will then asymptotically decay as a constant times $\exp (- \lambda t)$).

One useful additional remark is that one can couple two processes starting from $(0,0)$ and $(\eps, 0)$ using the same radial SLE$_6$. When $\eps$ is small, the probability that one of the two processes gets killed strictly before the other is bounded by the probability that three disjoint arms land exactly on the small boundary arc of length $\eps$, which is easily seen to be bounded by a constant times $\eps^2$ (this is one of the exponents that can be determined without SLE; see, e.g., \cite{Wlectures}). This shows that $G(\eps, 0) - G (0, 0) = O( \eps^2)$ as $\eps \to 0$, which corresponds to the ``balance'' condition mentioned above. For the other values of $\lambda$, the ``birth/killing'' rate at $(0,0)$ kicks in and this difference would be of order $\eps^{1/3}$.

\subsubsection*{The solutions and how to deduce the exponent}
AI suggests that the following functions (for all $\lambda \in (0,2/3)$) are positive eigenfunctions of the above operator with eigenvalue $- \lambda$, which satisfy the above boundary conditions up to a multiplicative constant:
\begin{equation}\label{explicitG}
{ G_\lambda(\alpha,\beta)}
 = {\sin(\alpha/2)^{1/3}}
 \int_\beta^{2\pi-\alpha}
 \frac{1}{[\sin(b/2)\sin((\alpha+b)/2)]^{2/3}} T_\lambda (\alpha, b) db
 \end {equation}
 for
 $$ T_\lambda (\alpha, b) :=
 \frac73\int_0^{\alpha/2+b}
 S_\omega(\alpha/2+b-u)(1-\cos u)^{2/3} du
 -\int_{\cos(\alpha/2)}^1
 (t-\cos(\alpha/2+b))^{2/3}f_{\omega}(t)dt,
 $$
where
$$
\omega:=\frac{\sqrt{1+12\lambda}}3,\quad S_\omega (u):= \frac {\sin (\omega u)}{\omega} , \quad
 f_{\omega}(t):=
{\null}_2F_1 \left(\tfrac53-\omega,\tfrac53+\omega;
                 \tfrac{13}{6};\tfrac{1-t}{2}\right)
.$$

Here ${}_2F_1$ denotes the Gauss hypergeometric function. With this normalization, $G_\lambda(0,\beta) = 2B(\tfrac13,\tfrac13)\, T_\lambda(0,2\pi)$ for all $\beta\in[0,2\pi)$, where $B$ is the beta function, so that $G_\lambda/G_\lambda(0,0)$ satisfies the normalization $G(0,\beta)=1$ above.

As argued above, one can then identify the value of $\lambda$ corresponding to the backbone exponent because of the condition
$G_\lambda (\alpha, 0) - G_\lambda (0,0) = o(\alpha^{1/3})$
as $\alpha \to 0$.
Indeed, it turns out that
$$
G_\lambda (\alpha, 0) = G_\lambda (0,0) +c(\lambda) \alpha^{1/3}
+o(\alpha^{1/3}),
$$
where $c(\lambda)$ can be computed explicitly and vanishes at only one value.

Let us briefly indicate how to derive this, since this discussion is quite instructive of the type of arguments and manipulations that come into play, and also explains the particular form of the backbone exponent:
Note first that $$T_\lambda(0,b)=j_\lambda(b):=\frac73\int_0^b S_\omega(b-u)(1-\cos u)^{2/3}\,du.$$
When $\beta=0$ and $\alpha\to0$, the $\alpha^{1/3}$ term comes from the prefactor $\sin(\alpha/2)^{1/3}$ in (\ref{explicitG}), so one has to focus on the integral. The remaining integrand converges to $\sin(b/2)^{-4/3}j_\lambda(b)$, which is not integrable at $b=2\pi$; the part of the integral where $2\pi-b$ is of order $\alpha$ produces $G_\lambda(0,0)$ together with a compensating term of order $\alpha^{1/3}$. Combining the two contributions shows that $c(\lambda)$ is the Hadamard finite part
$$
c(\lambda)=2^{-1/3}\lim_{\eps\to0}\Big(\int_0^{2\pi-\eps}\frac{j_\lambda(b)}{\sin(b/2)^{4/3}}\,db-\frac{3\cdot2^{4/3}\,j_\lambda(2\pi)}{\eps^{1/3}}\Big).
$$
Some little manipulations [inserting the definition of $j_\lambda$, exchanging the order of integration, using Euler's beta integral
$\int_0^{2\pi-u}\bigl(\sin(t/2)/\sin((t+u)/2)\bigr)^p\,dt=2\pi\sin(p(\pi-u/2))/\sin(\pi p)$
(it is valid for $p\in(-1,1)$) and then replacing the exponent $4/3$ by $p$ in both powers of the sin functions and using analytical continuation to $p=4/3$] lead to
$$
c(\lambda)=\frac{56\cdot2^{1/3}\,\pi}{3\sqrt3\,(4\lambda-1)}\Bigl(S_\omega(2\pi)+\frac{3\sqrt3}{4}\Bigr).
$$
Note that $4S_\omega(2\pi)+ {3\sqrt3}$ has in fact two zeroes in $(0,2/3)$, one of which being $\omega = 2/3$, i.e.,  $\lambda=1/4$ (which is the polychromatic two-arm exponent), and the other one corresponding to some value of $\lambda$ in $(1/4, 2/3)$. We can either rule out $\lambda = 1/4$ by the a priori results from \cite {BN}, or actually note that when $\lambda = 1/4$, the zero is cancelled by the factor $4\lambda-1$ in the expression of $c(\lambda)$, so that $c(1/4)$ is indeed not equal to $0$. Hence, for $\lambda\in(0,2/3)$, the corner condition $c(\lambda)=0$ holds only for the unique root in  $(1/4,2/3)$ -- which is therefore the backbone exponent.

Writing out the entire self-contained proof (starting from the arguments in \cite {LSW5} up to the details about solving the PDE and the identification of the backbone exponent) ends up being a matter of 10 to 15 pages or so (half of this for the detailed derivation of the PDE outlined above and the other half to check that the above function is indeed the solution).

\subsection {The (generalized) annulus crossing functions}
\label {S22}
Let us briefly explain how things work out for crossings of annuli.

\subsubsection*{The PDE obtained along the lines of \cite {Du}}
Let $u(m,\theta) = F( \exp (-m), \theta)$. Simply writing out the martingale property for the crossing event along the exploration via an SLE$_6$ leads to
$$
 \partial_m u=3 \partial_\theta^2 u +H(m,\theta) \partial_\theta u
$$
in $(0, \infty)$$ \times (0, 2 \pi)$, where
$$
 H(m,\theta)=\cot(\theta/2)+4\sum_{n=1}^{\infty}
                   \frac{\sin(n\theta)}{e^{2nm}-1}$$
(note that $H=2\partial_\theta\log\Theta$, where $\Theta$ is defined below).
                   The boundary conditions are $u ( m, 2 \pi-) = 0$, $u(0+, \theta) = 1$, and $u(m, \theta) - u(m,0) = O(\theta^2)$ as $\theta \to 0$ (this last estimate comes from the half-plane three-arm exponent).
                   The fact that for the solution $u$ to this PDE, $u(\log (1/r), 0)$ is the annulus crossing probability follows from the annular exploration framework of \cite {Du}.
\subsubsection*{The explicit solution}
The function $u$ on $(0,\infty)$$ \times (0, 2 \pi)$ that solves the PDE turns out to be
$$
u (m , \theta) = \frac{\pi\sqrt3}{2B(1/3,1/3)}E(m)^2\Theta(m,\theta)^{1/3}\times\int_0^\theta
 \frac{g(m,t-\theta/2)}
      {[\Theta(m,t)\Theta(m,\theta-t)]^{2/3}}dt,
$$
where
 $$ \Theta(m,\theta) = 2\sum_{n=0}^{\infty}(-1)^n
 e^{-m(n+1/2)^2}\sin((n+1/2)\theta) $$
 and
$$ g(m,z) =\frac2\pi\sum_{n=0}^\infty
      \cos\bigl((n+\tfrac12)z\bigr)
      \cos\bigl((n+\tfrac12)\pi/2\bigr)
      e^{-3(n+1/2)^2m/4}.$$
When $\theta \to 0$, one can indeed check (with enough proficiency with theta functions...) that
$$ u(m, 0) = \frac{\pi\sqrt3}{2E(m)}g(m,0)=F(e^{-m}),$$
where $F$ is given by Cardy's formula (\ref{Cardy}) (the first identity follows from $\int_0^1[v(1-v)]^{-2/3}\,dv=B(1/3,1/3)$ together with Jacobi's identity $\partial_\theta\Theta(m,0)=E(m)^3$, and the second one from the Jacobi triple product identity and the modular transformation $\eta(i/t)=\sqrt{t}\,\eta(it)$), so that the annulus crossing probability formula is a particular case of this general expression for the function $u(m, \theta)$.

\subsection*{Other examples}
In fact, similar explicit expressions can be obtained for functions with two arms (either monochromatic or bichromatic) and two boundary arcs leading, for instance, to annulus versions of the previous backbone probabilities (i.e., the probability that two disjoint open arms join the inner and outer boundaries); in the annulus backbone case, however, the identification of the boundary conditions that correctly encode the two-arm event in the annulus requires some additional care. In the end, this for instance provides the formulas for the corresponding multiple crossing probabilities of annuli with some marked boundary arcs that generalize CLE crossing formulas (with no boundary arc) derived in the recent CLE integrability papers such as \cite {integ2,Misc3}.

\section {Some concluding remarks}
We conclude with the following remarks:
\begin {enumerate}
\item
It is natural to expect that similar proofs exist for most other integrability results for CLE, and that the various expressions that show up in these direct proofs might be interpreted in LQG terms. This could actually also provide routes to SLE-based derivations of the landmark Liouville formulas, which would nicely close the circle and join the various approaches in a harmonious way.

Recall that it is known (based on the CLE connection probabilities determined in \cite {MW}, or on the CLE percolation formulas from \cite {MSW2} combined with \cite {MSW}) which CLE$_\kappa$ is the only possible conformally invariant scaling limit of which critical FK($q$) or $O(N)$ model (this is similar to Schramm's identification of SLE$_6$ as the only candidate for the scaling limit of percolation interfaces).
\item Exploring a CLE$_\kappa$ via SLE$_\kappa (\rho, \kappa-6- \rho)$-type curves (see, e.g., \cite{WW13,MSW}) that progressively discover the CLE loops one by one naturally leads to ``intermediate'' configurations with CLEs with three marked boundary points (the boundary being then divided into three parts: the left side of the currently discovered loop, the right side of the currently discovered loop, and the ``untouched'' boundary part), which, for instance, correspond to the two angles $\alpha$ and $\beta$ in the backbone PDE. So, in order to make this SLE martingale approach work, one needs to compute a function $\tilde H$ of $N+2$ variables if $N$ is the number of marked points that the function is following: instead of looking for a function $H(x_1, \ldots, x_N)$, one ends up having to find a function $\tilde H ( \alpha, \beta, x_1, \ldots, x_N)$. The general philosophy here is that whenever $H$ has an explicit expression, $\tilde H$ might have one as well, which would lead to some CLE-type proof. This goes, for instance, along the lines of \cite {BB,DubEuler}, but with very explicit expressions for the functions $\tilde H$ and CLEs instead of SLEs. We leave it to others to ``prompt'' this for the various existing cases mentioned above or for new ones.
\item It would be instructive to see how the central charge (related to the restriction properties from \cite {LSWr,ShW}) enters the game in these equations.
 \item One can reflect on the role played by AI in the current state of the art that we are describing here. Note that the ``thinking/conceptualization'' remained human-based, while AI is left with the task of solving a PDE.
 But to make it really useful for our human understanding, it is of course desirable to work out the ``route'' to obtaining these explicit formulas. In the above two examples, this appears to be possible, even if it requires some juggling with theta functions and various changes of variables (see for instance the little discussion to identify the backbone exponent). It would of course also be very interesting to compare this with the CFT approaches, for instance, to see if there are some structural parallels and/or if some of the functions that appear in the these integrals have a LQG interpretation.
 \item Even if solving these PDEs appears, a posteriori, to be something that could have been achieved by some humans (for instance, one can wonder what would have happened if Don Zagier had looked into this particular equation -- recall that Kleban and Zagier observed the appearance of modular forms in crossing probability formulas \cite {KZ}) one can also argue that finding/guessing these explicit solutions at the PDE level is difficult, and that, in some sense, the route toward integrability found via CFT-type arguments and structures is more amenable to (clever) humans!
 \item One can also feel fortunate that all these integrability results were first obtained by other means, which brought with them many new ideas and structures. Otherwise, chances are that the results (with proofs along the lines presented here) would have been part of the large batch of OpenAI ``new results'' from a couple of weeks ago.
  \item In retrospect, one can also wonder how differently this field would have evolved if this ``proficient'' person had shown up 25 years ago: what would have been gained, but also how much would have been lost! This is all food for thought about how to set up future AI-assisted research.
\end {enumerate}

\subsection*{AI disclosure.} The use of AI in the preparation of the math part of this paper was carefully described within the paper (i.e., it produced the explicit formulas for the solutions to the PDEs and helped justify that they are indeed solutions). We also used Claude for proofreading.

\subsection*{Acknowledgments.} J.M.\ received support from ERC Consolidator Grant ARPF (Horizon Europe UKRI G120614).  The research of W.W.\ has been supported by a Research Professorship of the Royal Society.

\printbibliography

\end{document}